\def\End{{\rm End}}

\newtheorem{prop}{Proposition}
\newtheorem{theorem}{Theorem}
\newtheorem{definition}{Definition}

\documentstyle[12pt]{article}

\begin{document}
\title{The moduli space of complex Lagrangian submanifolds}
\author{N.J.HITCHIN\\[5pt]
\itshape  Mathematical Institute\\
\itshape 24-29 St Giles\\
\itshape Oxford OX1 3LB\\
\itshape England\\[12pt]}
\maketitle
\bigskip

\section{Introduction}
Developments in string theory over the past few years (e.g.\cite{SW},\cite{D}) have focussed attention on a differential geometric structure induced on the base space of an algebraically  completely integrable Hamiltonian system. This has been recently formalised by D.Freed \cite {F} as a {\it special K\"ahler structure}.

The purpose of this paper is to provide an alternative approach to the geometry of special K\"ahler manifolds, one that is motivated by the desire to understand a more general situation than that afforded by integrable systems. We seek the natural geometrical structure on the moduli space $M$ of deformations of a compact, complex Lagrangian submanifold $Y$ in a complex K\"ahlerian symplectic manifold $X$. In many respects what we do parallels the approach of an earlier paper \cite{H} which began an investigation into the geometry of the moduli space of compact special Lagrangian submanifolds in a Calabi-Yau manifold. We are essentially attacking the special case where the Calabi-Yau is hyperk\"ahler, though we shall not need the full force of the existence of a hyperk\"ahler metric on $X$.

Our viewpoint, as in \cite{H}, is to pay less attention to  the holomorphic structure of the situation, and more to the symplectic one. Thus a complex Lagrangian submanifold  of $X$ can be characterized as a real submanifold on which the real and imaginary parts of the holomorphic symplectic 2-form vanish. Correspondingly, the differential geometric structure on the moduli space $M$ is induced from a local embedding of $M$ into $H^1(Y,{\bf R})\times H^1(Y,{\bf R})$ which is Lagrangian with respect to two natural constant symplectic forms. We show that this ``bilagrangian'' condition for a submanifold of a product $V\times V$ of real symplectic vector spaces is equivalent to the structure of a special (pseudo-) K\"ahler metric on $M$. Moreover, it is easy to see from this point of view that a choice of symplectic basis of $V$ yields the known fact that any special K\"ahler metric is generated by a single holomorphic function -- the holomorphic prepotential.

Finally, we derive from our formalism the hyperk\"ahler metric introduced in a string-theoretic context several years ago by Cecotti, Ferrara and Girardello \cite {CFG}. It can be seen as a special case of the Legendre transform construction of Lindstr\"om and Ro\v cek \cite {LR}, and yields, in the context of our moduli space,  a hyperk\"ahler metric on an open set of Markman's moduli space of Lagrangian sheaves \cite{Mark}.

Our approach offers a different perspective to special K\"ahler geometry, and in particular draws attention to a single naturally defined function $\phi$ which plays an important role: it is the Hamiltonian for the  fundamental vector field, a potential for the K\"ahler metric and, with respect to one of the complex structures, a potential also for the associated hyperk\"ahler metric.
\vskip .25cm
The author wishes to thank Dan Freed for introducing the subject to him and the Institute for Advanced Study for its hospitality.

\section{Complex Lagrangian submanifolds}
Let $X$ be a complex symplectic manifold of complex dimension $2n$. It has a holomorphic symplectic 2-form $\omega^c$ which we write in terms of its real and imaginary parts:
$$\omega^c=\omega_1+i\omega_2$$
These two closed forms are real symplectic forms and define the structure of a complex symplectic manifold on $X$. We see this as follows. Given a closed form $\omega^c$, we consider the distribution in $E\subset T\otimes{\bf C}$ defined by the complex vector fields $U$ with $\iota(U)\omega^c=0$. If $\omega^c$ satisfies the algebraic condition that $E\oplus\bar E=T\otimes{\bf C}$ then it defines an almost complex structure.
 This is integrable because if $\iota(U)\omega^c=0$,
$${\cal L}_U\omega^c=d(\iota(U)\omega^c)+\iota(U)d\omega^c=0$$
so if $U$ and $V$ are sections of $E$, ${\cal L}_U\omega^c=0$ and $\iota(V)\omega^c=0$, so
$$\iota([U,V])\omega^c={\cal L}_U(\iota(V)\omega^c)-\iota(V)({\cal L}_U\omega^c)=0$$  
and hence $[U,V]$ is a section of $E$. 

Similarly we have the following
\begin{prop} A real $2n$-dimensional submanifold $Y\subset X$ is complex Lagrangian if and only if $\omega^c\vert_Y=0$.
\end{prop}
\noindent {\it Proof:} If the submanifold $Y$ is complex Lagrangian, then $\omega^c\vert_Y=0$ by definition. Conversely, if $\omega^c\vert_Y=0$, we need to show that $Y$ is a complex submanifold, that its tangent spaces are complex.
Now the  complex structure $I$ on $X$ is defined algebraically by the two real symplectic forms $\omega_1,\omega_2$. Instead of the description above,  we can think of   each giving an isomorphism $\varphi_i:TX\cong T^*X$ and then $I=\varphi_2^{-1}\varphi_1$. Since $Y$ is Lagrangian with respect to both symplectic forms, then both $\varphi_1$ and $\varphi_2$ map $TY$ isomorphically to the conormal bundle $N^*Y\subset T^*X$. thus $I=\varphi_2^{-1}\varphi_1$ preserves $TY$. Hence $Y$ is a complex submanifold, $\omega^c\vert_Y=0$ and so $Y$ is complex Lagrangian.
\vskip .25cm
We lose nothing therefore by focussing our attention on the symplectic aspects of complex Lagrangian submanifolds. We do, however, need to know that there is a good local moduli space of deformations of $Y\subset X$. In general  deformations of compact complex submanifolds can be obstructed, but it follows from the paper of Voisin \cite{V} that this is not the case when the submanifold is Lagrangian. In fact if $X$ has a  hyperk\"ahler metric, this is also a consequence of the differential geometric argument of McLean \cite{Mac}. There thus exists a local moduli space $M$, which is a complex manifold, and such that there is a natural isomorphism
$$T_{[Y]}M\cong H^0(Y,N)$$
from the tangent space at the point $[Y]\in M$ representing $Y$ and the space of holomorphic sections of the normal bundle $N$ to $Y\subset X$. As $Y$ is complex Lagrangian, $\omega^c$ defines a holomorphic isomorphism from $N$ to $T^*Y$, so 
$$T_{[Y]}M\cong H^0(Y,T^*)$$
From the holomorphic point of view this is the tangent space to the moduli of deformations of $Y$ as a complex submanifold. The infinitesimal deformations as a complex {\it Lagrangian} submanifold correspond to those sections of the normal bundle for which the corresponding 1-forms in $H^0(Y,T^*)$  are {\it closed}. But
if we assume that $X$ has a K\"ahler form $h$ so that $Y$ is also K\"ahler, then any holomorphic 1-form is closed. Moreover  the real dimension of $M$ is then given by $b_1(Y)=2\dim_{\bf C} H^0(Y,T^*)$. We make this K\'ahlerian assumption from now on. 
We shall investigate next the local differential geometry of $M$.

\section {The moduli space}
Let $Z$ be a local universal family of deformations of the complex Lagrangian submanifold $Y_0\subset X$, so that $Z$ is a complex manifold with a holomorphic projection $\pi:Z\rightarrow M$ and a holomorphic map $F:Z\rightarrow X$ such that
$$F(\pi^{-1}([Y]))=Y$$
Consider the 2-form $F^*\omega_1$ on $Z$. If $x_1,\dots,x_{2n}, y_1,\dots,y_{2m}$ are real local coordinates on $Z$ with $y_1,\dots,y_{2m}$ coordinates on $M$ and $\pi(x_1,\dots,y_{2m})=(y_1,\dots,y_{2m})$ then since each fibre $Y$ is Lagrangian with respect to $\omega_1$, $F^*\omega_1\vert_Y=0$ and so
$$F^*\omega_1=\sum a_{ij}dx_i\wedge dy_j+\sum b_{ij}dy_i\wedge dy_j$$
Furthermore, since $F^*\omega_1$ is closed,
\begin{equation}
\sum \frac{\partial a_{ij}}{\partial x_k}dx_k\wedge dx_i=0
\label {d}
\end{equation}
We can see in concrete coordinate terms here that, for each $j$, the 1-form $\sum a_{ij}dx_i$ on $Y$ is closed. More invariantly, it says that if $U$ is a tangent vector to $M$ at $[Y]$, then if $\tilde U$ is a lift to a vector field along $Y$, the 1-form $(\iota(\tilde U)F^*\omega_1)\vert_Y$ is closed and independent of the choice of lifting.

From (\ref{d}), integrating $F^*\omega_1$ over two homologous 1-cycles in a fibre of $\pi$ gives the same result. Now working locally in $M$, we assume that $M$ is contractible, and so we can by homotopy invariance identify the homology of each fibre. Take a homology class $A\in H_1(Z,{\bf Z})\cong H_1(Y,{\bf Z})$ and choose a circle fibration in $\pi:Z\rightarrow M$ such that each fibre is in the class $A$. Integrating along the fibres, we obtain canonically a closed 1-form
$$(\pi_A)_*F^*\omega_1\in \Omega^1(M)$$
and since $M$ is assumed contractible, a smooth function $\mu_A$ on $M$, well-defined up to an additive constant, such that $d\mu_A=(\pi_A)_*F^*\omega_1$. Putting all the functions together gives a map
$$\mu:M\rightarrow H^1(Y,{\bf R})$$
This function by definition has the property that $d\mu (U)$ is the cohomology class of the closed form $(\iota(\tilde U)F^*\omega_1)\vert_Y$.

\begin{prop} $\mu$ is a local diffeomorphism.
\end{prop}
\noindent {\it Proof:} We think in  terms of the holomorphic fibration $\pi :Z\rightarrow M$ and a tangent vector $U$ at $[Y]$. The  1-form $(\iota(\tilde U)F^*\omega^c) \vert_Y$ on $Y$ is independent of the choice of lift and thus, taking local holomorphic lifts, is a well-defined global holomorphic 1-form. This is the canonical isomorphism $T_{[Y]}M\cong H^0(Y,T^*)$. Now if $d\mu(U)=0$, by the definition of $\mu$, the cohomology class of the real part of $\theta=(\iota(\tilde U)F^*\omega^c) \vert_Y$ is zero. But $Y$ is  a K\"ahler manifold so we have $H^1(Y,{\bf C})\cong H^{1,0}\oplus H^{0,1}$ and $\theta +\bar \theta$ cannot be cohomologically trivial unless $\theta=0$. This means that $(\iota(\tilde U)F^*\omega^c)\vert_Y=0$ and $U=0$. Thus by the inverse function theorem, $\mu$ is a local diffeomorphism.
\vskip .25cm
Similarly, using the other symplectic form $\omega_2$, we get a map
$$\nu:M\rightarrow H^1(Y,{\bf R})$$
and, put together, a smooth map
$$w=(\mu,\nu):M\rightarrow H^1(Y,{\bf R})\times H^1(Y,{\bf R})$$
Thus $w(M)$ is a smooth submanifold such that the projection onto each factor is a local diffeomorphism.
 \vskip .25cm
 The vector space $H^1(Y,{\bf R})$ has a real constant symplectic form defined by the restriction of the K\"ahler form $h$ on $X$:
$$\langle a,b\rangle=\int_Y \alpha\wedge \beta \wedge h^{n-1}$$
where $\alpha,\beta$ are representative 1-forms for $a$ and $b$. This clearly only requires the cohomology class of the K\"ahler form $h$.
We  define two constant symplectic forms on $H^1(Y,{\bf R})\times H^1(Y,{\bf R})$:
\begin{eqnarray}
\Omega_1((a_1,a_2),(b_1,b_2))&=&\langle a_1,b_2\rangle+\langle a_2,b_1\rangle
\label {o1}\\
\Omega_2((a_1,a_2),(b_1,b_2))&=&\langle a_1,b_1\rangle-\langle a_2,b_2\rangle 
\end{eqnarray} \label {o2}
Take a basis for $V$, so that the skew form has matrix $\omega_{ij}$, then in the corresponding linear coordinates
\begin{eqnarray}
\Omega_1&=&2\sum \omega_{ij}dx_i\wedge dy_j \label{o3}\\
\Omega_2&=& \sum \omega_{ij}dx_i\wedge dx_j-\sum \omega_{ij}dy_i\wedge dy_j \label{o4}
\end{eqnarray}
We now have
\begin{theorem} $w(M)\subset H^1(Y,{\bf R})\times H^1(Y,{\bf R})$ is Lagrangian with respect to $\Omega_1$ and $\Omega_2$.
\end{theorem}
\noindent {\it Proof:} The holomorphic symplectic form $\omega^c$ is of type $(2,0)$ so since $F$ is holomorphic, $F^*(\omega^c)^2$ has type $(4,0)$. If $U$ is a tangent vector on $M$ at $[Y]$, then as we have seen, $(\iota(\tilde U)F^*\omega^c )\vert_Y$ is independent of the choice of lifting $\tilde U$, because $\omega^c \vert_Y=0$. Now
$$\iota(\tilde U)\iota(\tilde V)F^*(\omega^c)^2=2(\iota(\tilde U)\iota(\tilde V)F^*\omega^c) (F^*\omega^c) +2(\iota(\tilde U)F^*\omega^c )(\iota(\tilde V)F^*\omega^c)$$
and restricting to $Y$, 
\begin{equation}
\iota(\tilde U)\iota(\tilde V)F^*(\omega^c)^2 \vert_Y=2(\iota(\tilde U)F^*\omega^c )\vert_Y \wedge  (\iota(\tilde V)F^*\omega^c) \vert_Y
\label{20}
\end{equation}
But the left hand side is of type $(2,0)$, hence
$$\iota(\tilde U)\iota(\tilde V)F^*(\omega^c)^2 \vert_Y \wedge h^{n-1}\vert_Y$$
is of type $(n+1,n-1)$ and so vanishes since $Y$ is complex of dimension $n$. Hence from (\ref{20})
$$\int_Y \iota(\tilde U)F^*\omega^c \wedge \iota(\tilde V)F^*\omega^c \wedge h^{n-1}=0$$
but this means that
$$\langle (d\mu+id\nu)(U),(d\mu+id\nu)(V)\rangle=0$$
and so equating to zero real and imaginary parts,
\begin{eqnarray*}
\langle d\mu(U),d\mu(V)\rangle-\langle d\nu(U),d\nu(V)\rangle&=&0\\
\langle d\nu(U),d\mu(V)\rangle+\langle d\mu(U),d\nu(V)\rangle&=&0
\end{eqnarray*}
These two conditions are precisely the vanishing of $\Omega_2$ and $\Omega_1$ respectively on $w(M)$.
\vskip .25cm
\noindent {\bf Remark:} Theorem 1 demonstrates that the structure of the moduli space -- defined as a submanifold on which two symplectic forms vanish -- parallels the structure of the objects it parametrizes. This is also the philosophy behind the description in \cite{H} of the moduli space of special Lagrangian submanifolds of a Calabi-Yau manifold.

\section{Special K\"ahler manifolds}
We shall show that, as a consequence of the ``bilagrangian'' property of Theorem 1, $M$ inherits a special K\"ahler structure. The definition, as given in \cite{F}, is the following:
\begin{definition} A special K\"ahler manifold is a complex manifold $M$ with
\begin{itemize}
\item
a K\"ahler metric $g$ with K\"ahler form $\omega$
\item
a flat torsion-free connection $\nabla$ such that 
\item
$\nabla \omega=0$ and
\item
$d_{\nabla}I=0\in \Omega^2(M,T)$
\end{itemize}
\end{definition}
To clarify the last property, we think of the complex structure $I$, an endomorphism of the tangent bundle $T$, as a 1-form with values in $T$, i.e. $I\in \Omega^1(M,T)$. The connection $\nabla$  defines a covariant exterior derivative $d_{\nabla}:\Omega^p(M,T)\rightarrow \Omega^{p+1}(M,T)$ and we require $d_{\nabla}I=0$. This is weaker than $\nabla I=0$ and indeed, since $\omega$ and $I$ determine $g$, the latter condition would say that $\nabla$ is the Levi-Civita connection, and all we would be looking at is a flat K\"ahler manifold. The reader should be warned that special K\"ahler manifolds do not form a very interesting class of global Riemannian structures -- it has been shown by Lu \cite{Lu} that any {\it complete} special K\"ahler manifold is flat.
\vskip .25cm
Now let $V$ be a real symplectic vector space with skew form $\langle\,,\rangle$. As in (\ref{o1}),(\ref{o2}) we define two constant symplectic forms $\Omega_1,\Omega_2$ on $V\times V$. We can also define an indefinite metric $g$ on $V\times V$ by
$$g((a,b),(a,b))=\frac{1}{2}\langle a,b\rangle$$
Then we have the following

\begin{theorem}  Let $M\subset V \times V$ be a submanifold which is Lagrangian
for $\Omega_1$ and $\Omega_2$, and transversal to the two projections onto $V$. Then
$g\vert_M$ is a special pseudo-K\"ahler metric. Conversely, any special pseudo-K\"ahler metric on a manifold $M$ is  locally induced from an 
embedding in $V\times V$ for some $V$.
\end{theorem}

\noindent{\it Proof:} It will be convenient to use the nondegenerate form $\langle \, ,\rangle$ on $V$ to identify $V$ with $V^*$. Under this identification, $\Omega_1$ becomes essentially the canonical symplectic form on $T^*V\cong V\times V^*$, by setting $d\xi_i=\sum \omega_{ij}dy_j$ in (\ref{o3}),(\ref{o4}). We then have the following expressions for $\Omega_1$ and $\Omega_2$:
\begin{eqnarray}
\Omega_1&=&2\sum dx_i\wedge d\xi_i \label {O1}\\
\Omega_2&=& \sum \omega_{ij}dx_i\wedge dx_j+\sum \omega^{ij}d\xi_i\wedge d\xi_j \label{O2}
\end{eqnarray}
\vskip .25cm
To begin, we use the projection onto the first factor $V$ to locally identify $M$ with a flat symplectic vector space. This provides us with our flat connection $\nabla$ with $\nabla \omega=0$.  If we use the coordinates $x_1,\dots,x_{2n}$, then covariant derivatives using $\nabla$ are just ordinary derivatives.

Now since $M$ is Lagrangian with respect to the canonical symplectic form $\Omega_1$ on $T^*V$
and transversal to the projection to $V$,  the embedding is defined by the graph of the derivative of a function on $V$ so in coordinates 
$$\xi_j=\frac{\partial \phi}{\partial x_j}$$
for some function $\phi(x_1,\dots,x_{2n})$. The tangent vector $\partial/\partial x_j$ of
$M$ then lies in $V\times V^*$ as the vector
$$X_j=\frac{\partial}{\partial x_j}+\sum \frac{\partial \xi_k}{\partial
x_j}\frac{ \partial}{\partial \xi_k}=\frac{\partial}{\partial x_j}+\sum_k
\frac{\partial^2 \phi}{\partial x_k \partial x_j}\frac{ \partial}{\partial
\xi_k}$$
Since the metric on $V\times V^*$  is defined by 
$$g((x,\xi),(x,\xi))=\frac{1}{2} \langle x,\xi \rangle$$
 the induced metric on $M$ is 
$$\sum g_{kj}dx_kdx_j=\sum g(X_k,X_j)dx_kdx_j=\sum \frac{\partial^2 \phi}{\partial x_k
\partial x_j}dx_kdx_j$$
In general this metric may not be  positive definite. It is nondegenerate however, for if $\sum g_{ij}a_j=0$, then $\sum a_jX_j=\sum a_j\partial/\partial x_j$ so that projection of this tangent vector onto the second factor in $V\times V$ is zero. By the transversality assumption this means each $a_j=0$. 
\vskip .25cm
Consider now the second Lagrangian condition: $\Omega_2$
vanishes on $M$. From (\ref{O2}) this says, using $\xi_k=\partial \phi/\partial x_k$ and $g_{kj}={\partial^2
\phi}/{\partial x_k \partial x_j}$,
$$0=\Omega_2(X_k,X_j)=\omega_{kj}+\sum g_{ka}g_{jb}\omega^{ab}$$
or, writing $I^j_k=\sum\omega^{ja}g_{ak}$, that $I^2=-1$. This is the almost complex structure. Since
$\omega_{kj}=-\omega_{jk}$, $I$ is skew adjoint with
respect to $g$. Now let $X$ be the Hamiltonian vector field for the function $\phi$.
 We have 
$$X=\sum a_j{\frac{\partial }{\partial x_j}}= \sum \omega^{ij}{\frac{\partial \phi}{\partial x_i}}{\frac{\partial }{\partial x_j}}$$
so
$$\frac{\partial a_j}{\partial x_k}=\sum \omega^{ij}g_{ik}=I^j_k$$
Hence $I=d_{\nabla}X$ and so $d_{\nabla}I=d_{\nabla}^2X=0$ since the connection is flat. Hence $I$ satisfies the compatibility condition with the flat connection.
\vskip .25cm
\noindent It remains to show that $I$ is integrable. But consider the complex
functions
\begin{equation}
z_j=x_j-i\sum \omega^{jk}\frac{\partial \phi}{\partial x_k}
\label{forms}
\end{equation}
Differentiating, we obtain the $2n$ complex 1-forms
$$dz_j=dx_j-i\sum\omega^{jk}\frac{\partial^2 \phi}{\partial x_k \partial
x_l}dx_l=dx_j-i\sum \omega^{jk}g_{kl}dx_l=dx_j-i\sum I^j_l dx_l$$
and these are clearly of type $(1,0)$. We need to find $n$ linearly independent ones.  Let $E\subseteq \Lambda^{1,0}$ be the distribution spanned by $dz_1,\dots,dz_{2n}$. Then  for each $j$, $2dx_j=dz_j+d\bar z_j$, and since $dx_1,\dots,dx_{2n}$ forms a basis, $E\oplus \bar E=\Lambda^1$ and the rank of $E$ is $n$. 

Thus the metric $g$ and the connection $\nabla$ satisfy all the conditions for a special pseudo-K\"ahler (i.e. possibly indefinite) metric.
\vskip .25cm
Now consider the converse.  Let $M$ be a special pseudo-K\"ahler manifold, and $(x_1,\dots,x_{2n})$ flat
local coordinates, so that the covariant derivative is the ordinary derivative and the coefficients of the symplectic form are constant. Consider the 1-form 
$$\alpha_k=\sum \omega_{kl}I^l_jdx_j$$
Now 
$$d\alpha_k=\sum \omega_{kl} \frac{\partial I^l_j}{\partial x_m}dx_m\wedge dx_j=0$$
since $\omega_{kl}$ is constant and $d_{\nabla}I=0$. Thus locally there are
functions $\xi_k$ such that
$$\alpha_k=d\xi_k.$$ 
We map $M$ to ${\bf R}^{2n}\times {\bf R}^{2n}$ by 
$$(x_1,\dots,x_{2n})\mapsto (x_1,\dots,x_{2n},\xi_1,\dots,\xi_{2n})$$
\vskip .25cm
\noindent First we claim the image  is Lagrangian with respect to the symplectic
form $\Omega_1=2\sum dx_j\wedge d\xi_j$.  But  restricting $\sum dx_j\wedge d\xi_j$ gives 
$$
dx_j\wedge \alpha_j=\sum dx_j\wedge \omega_{jl}I^l_kdx_k=\sum dx_j\wedge g_{jk}
dx_k=0$$
since the metric tensor $g_{jk}$ is symmetric.
\vskip .25cm
\noindent Next consider the symplectic form $\Omega_2 =\sum
\omega_{ij}dx_i\wedge dx_j+\omega^{ij}d\xi_i\wedge d\xi_j$. Restricted to $M$
this is 
$$ \sum \omega_{ij}dx_i\wedge dx_j+\omega^{ab}g_{ai}dx_i\wedge g_{bj}dx_j$$
But since $I^2=-1$ this too is zero. 
\vskip .25cm
Hence taking $V$ to be ${\bf R}^{2n}$ with the skew form $\omega_{ij}$, we obtain
a  local embedding of $M$ in $V\times V$, Lagrangian with respect to both forms. It is straightforward to see that the induced metric is $g_{ij}$.
\vskip .5cm
\noindent One of the features of the above approach is the fundamental role of the function $\phi$. Here is another aspect of this:
\begin{prop} The function $\phi$ is a K\"ahler potential.
\end{prop}

\noindent{\it Proof}: Consider 
$$d(Id\phi)=\sum \frac{\partial }{\partial x_k}(I^i_j \frac{\partial \phi}{\partial x_i}) dx_k\wedge dx_j$$
Now 
$$\frac{\partial }{\partial x_k}I^i_j=\sum \frac{\partial }{\partial x_k}(\omega^{ia} g_{aj})=\sum \omega^{ia} \frac{\partial^3 \phi}{\partial x_a \partial x_j\partial x_k} $$
so
$$d(Id\phi)=\sum \omega^{ia} \frac{\partial^3 \phi}{\partial x_a \partial x_j\partial x_k}\frac{\partial \phi}{\partial x_i} dx_k\wedge dx_j+I^i_j \frac{\partial^2 \phi}{\partial x_i\partial x_k} dx_k\wedge dx_j$$
The first term vanishes by the symmetry in $j,k$ and the second term is
$$\sum I^i_j g_{ik} dx_k\wedge dx_j=-\sum \omega_{kj}dx_k\wedge dx_j=-\omega$$
where $\omega$ is the K\"ahler form.
\vskip .5cm
\noindent We can now prove
\begin{theorem} The moduli space $M$ of deformations of a complex Lagrangian submanifold of a complex K\"ahlerian symplectic manifold $X$ has a naturally induced special K\"ahler structure.
\end{theorem}
\noindent{\it Proof:} We have to show firstly that $w(M)\subset H^1(Y,{\bf R})\times H^1(Y,{\bf R})$ satisfies the transversality of the previous theorem, but this follows from Proposition 2 applied to $\mu$ and $\nu$. We also need the metric to be  definite. But
$$g(X,X)=\langle d\mu(X),d\nu(X)\rangle=-\frac{i}{2}\langle (d\mu+id\nu)(X),(d\mu-id\nu)(X)\rangle$$
$$=-\frac{i}{2}\int_Y\iota(X)\omega^c\wedge \iota(X)\bar \omega^c\wedge h^{n-1}$$
which is definite. We should also show that the complex structure $I$ coincides with the natural complex structure of the moduli space when considered as a moduli space of complex submanifolds. From (\ref{forms}), the functions
$$z_j=x_j-i\sum \omega^{jk}\frac{\partial \phi}{\partial x_k}$$
are antiholomorphic with respect to $I$. But  $\xi_j=\partial \phi/\partial x_j$, so $z_j=x_j-i\sum \omega^{jk}\xi_k$. If we return from $V\times V^*$ to $V\times V$ using $\omega$ we see that this   function is obtained by taking a class $A_j\in H_1(Y,{\bf R})$ and forming
$$z_j=\langle \mu-i\nu, A_j\rangle$$
Since $d(\mu+i\nu)(U)=\iota(\tilde U)\omega^c$ is of type $(1,0)$, $dz_j=\langle d\mu-id\nu, A_j\rangle$  is of type $(0,1)$ with respect to the complex structure of the moduli space of compact complex submanifolds, so the two complex structures coincide.

Finally we should pass from the local to the global point of view on $M$, a point which is relevant in particular to the situation where $M$ is the base space of a completely integrable system. There is an additive ambiguity in the choice of the function $\mu$ but $d\mu$ gives an isomorphism between $TM$ and the trivial bundle $M\times H^1(Y,{\bf R})$, and this is the flat connection $\nabla$ of the special K\"ahler structure. Globally on $M$, the cohomology of the fibres of $Z\rightarrow M$ defines a flat vector bundle on $M$ (homotopy invariance of cohomology defines the Gauss-Manin connection) and so the isomorphism $d\mu$ provides a flat connection on $TM$. The symplectic form is preserved by the Gauss-Manin connection, and since the complex structure on $M$ is globally defined, so is the metric $g$, which is defined by $I$ and $\omega$.
\vskip .5cm
\noindent {\bf Remarks:}

\noindent 1. One of the well-known features of special K\"ahler geometry is the fact that any special K\"ahler metric is derived from a single holomorphic function ${\cal F}$ of $n$ variables. It is known as the holomorphic prepotential on $M$. This fact is rather easily seen using our bilagrangian formulism. For this purpose we choose a symplectic basis on $V$. The corresponding coordinates
 $x_1,\dots,x_{2n}$ give
$$\omega=\sum dx_j\wedge dx_{n+j}$$
and so 
\begin{eqnarray*}
\Omega_1&=&2\sum_1^{2n} dx_j\wedge d\xi_j \\
\Omega_2&=& \sum_1^n dx_j\wedge dx_{n+j}-\sum_1^n d\xi_j\wedge d\xi_{n+j} \
\end{eqnarray*}
and then
$$\Omega^c=\frac{1}{2}\Omega_1+i\Omega_2=\sum_i^n d(x_j+i\xi_{n+j})\wedge d(\xi_{j}+ix_{n+j})=\sum_1^n dv_j\wedge dw_{j}$$
We see that $\Omega^c$ can be identified with the canonical complex symplectic form on $T^*{\bf C}^n$. From Proposition 1, a submanifold on which $\Omega_1$ and $\Omega_2$ vanish is the same thing as a complex Lagrangian submanifold of $T^*{\bf C}^n$, but this is given by the graph of the derivative of a holomorphic function
\begin{equation}
v_{j}=\frac{\partial {\cal F}}{\partial w_{j}}
\label{FF}
\end{equation}
 From Theorem 2 this is all we need for a special K\"ahler manifold. 
 \vskip .15cm
 \noindent 2. In the standard presentation of the prepotential, its second derivative gives a holomorphic map from the base space of an integrable system to the Siegel upper half-space -- the moduli space of polarized abelian varieties, expressed as symmetric matrices with positive definite imaginary part. Such  a description involves choosing a symplectic basis for $H_1(Y,{\bf Z})$ (the classical $A$ and $B$ cycles) which is what we have done to introduce the holomorphic function ${\cal F}$. By contrast the real function $\phi$ requires no such choice. All we have chosen is projection onto the first factor in $V\times V$ to define $\phi$. We postpone the discussion of the relationship between $\phi$ and ${\cal F}$ to the next section, where we study some associated hyperk\"ahler constructions.
\vskip .15cm
\noindent 3. In the bilagrangian picture of $M\subset V\times V$ we get another flat torsion-free connection  by projecting onto the second factor. From the second Lagrangian condition 
$$0=\Omega_2\vert_M= \sum \omega_{ij}dx_i\wedge dx_j-\sum \omega_{ij}dy_i\wedge dy_j$$
 the pull back of the flat symplectic form under this projection is the same $\omega$. The function $\phi$ is then replaced by its Legendre transform and the new flat coordinates $dy_j$ are related to the old ones by $dy_j=Idx_j$. From the point of view of the moduli space of complex Lagrangian submanifolds, all we have done here is to replace $\omega_1$ by $\omega_2$, or what is essentially the same, to replace the complex symplectic form $\omega^c$ by $i\omega^c$. Clearly we can multiply $\omega^c$ by $e^{i\theta}$ and still have the same moduli space but now a family of flat connections parametrized by the circle. This is one viewpoint to the study of Higgs bundles as in Simpson's higher dimensional approach \cite{S}. In fact  any special  K\"ahler manifold can be thought of as having a Higgs bundle structure on $T\oplus T^*$. Recall that a Higgs bundle corresponding to a local system on a K\"ahler manifold $M$ consists of a holomorphic vector bundle $E$ with a unitary connection $A$ and a section $\Phi\in H^0(M, \End E\otimes T^*)$, the Higgs field, such that the connection
$$d_A+e^{i\theta}\Phi+ e^{-i\theta}\Phi^*$$
is flat for all $\theta$. In the case of a special K\"ahler manifold, the connection $A$ is the Levi-Civita connection on $T\oplus T^*$ and the Higgs field has the form
$$\pmatrix{0&0\cr
             \Theta &0\cr}$$
where $\Theta \in H^0(M,Sym^3 T^*)$ is a holomorphic cubic form. Since Simpson's original approach to Higgs bundles was derived from  variations  of Hodge structure, this  begins to take us back to the picture of a moduli space of complex manifolds which was the original motivation for this paper.
\section{Hyperk\"ahler metrics}

In \cite{CFG} Cecotti, Ferrara and Girardello showed how to define naturally a hyperk\"ahler metric on a certain bundle over a special K\"ahler manifold (see also \cite{F}). As we have seen a special K\"ahler metric can be defined via a single holomorphic function, so we have  a straightforward way of constructing hyperk\"ahler metrics. We shall show here that this construction is in fact a special case of an earlier technique called the {\it Legendre transform} construction of Lindstr\"om and Ro\v cek \cite{LR}, \cite{HKLR}.

Recall (see for example \cite{AH}) that a hyperk\"ahler metric is determined by three symplectic forms $\omega_1,\omega_2,\omega_3$ satisfying some algebraic conditions, namely that if $\varphi_i:T\rightarrow T^*$ is the isomorphism  determined by $\omega_i$, and we define $J_1=\varphi_3^{-1}\varphi_2,J_2=\varphi_1^{-1}\varphi_3, J_3=\varphi_2^{-1}\varphi_1$, then $J_1,J_2,J_3$ obey the quaternionic identities $J_1^2=J_2^2=J_3^2=-1$. We work locally first and let $M$ be a special K\"ahler manifold. On the product
$M\times {\bf R}^{2k}$
take the three  symplectic 2-forms $\omega_1,\omega_2,\omega_3$ defined by 
\begin{equation}
\omega_1=\sum \frac{\partial ^2 \phi}{\partial x_j \partial x_k}dx_j\wedge dy_k
\label{omega1}
\end{equation}
\begin{equation}
\omega_2+i\omega_3=-\frac{1}{2}\sum \omega_{jk}d(x_j+iy_j)\wedge d(x_k+iy_k)
\label{omegac}
\end{equation}
Then
$$J_3(\frac{\partial}{\partial x_i})=\sum
\omega^{jk}g_{ik}\frac{\partial}{\partial y_j}=\sum 
I^j_i\frac{\partial}{\partial y_j}$$
and
$$J_3(\frac{\partial}{\partial y_i})=\sum
\omega^{jk}g_{ik}\frac{\partial}{\partial x_j}=\sum 
I^j_i\frac{\partial}{\partial x_j}$$
 So in block matrix form
$$J_3=\pmatrix{0&I\cr
            I&0\cr}$$
Similarly
$$J_2=\pmatrix{-I&0\cr
             0&I\cr}$$
so that $J_3^2=J_2^2=-1$ and $J_2J_3=-J_3J_2$. Thus $J_2,J_3$ generate an
action of the quaternions, with $J_1=J_2J_3$ given by
$$J_1=\pmatrix{0&1\cr
             -1&0\cr}$$
 and so the symplectic forms $\omega_1,\omega_2,\omega_3$ define a hyperk\"ahler metric. 
\begin{prop} $\phi$ is a K\"ahler potential with respect to the complex structure $J_1$.
\end{prop}
\noindent{\it Proof:} From (\ref{omegac}) and the discussion in section 2, it is clear that $z_j=x_i+iy_j$ for $1\le j \le 2n$ are complex coordinates in this complex structure. But $\phi$ is independent of $y_j$ so 
$$\partial \bar \partial \phi= \sum \frac{\partial ^2 \phi}{\partial x_j \partial x_k}dz_j\wedge d\bar z_k=-2i\omega_1$$
from (\ref{omega1}).

\vskip .15cm
\noindent {\bf Remark:} Note that the projection from $M\times {\bf R}^{2n}$ to $M$ is holomorphic  in the complex structure $J_2$. Using the $\partial, \bar \partial$ operators in  that structure $\partial \bar \partial \phi$ gives the pull-back of the K\"ahler form on $M$, which is degenerate.  
\vskip .15cm

Since $\omega_{jk}$ is
constant and $\phi$ is independent of $y_j$, each symplectic form is invariant by the vector field $\partial/\partial
y_j$. We thus have a triholomorphic action of ${\bf R}^{2n}$ on $M\times {\bf R}^{2n}$.
The Legendre transform method is a canonical construction of hyperk\"ahler manifolds $X^{4n}$ with a triholomorphic action of ${\bf R}^n$, so in our case we have far more symmetry. Strictly speaking (a point not emphasized in the literature on this method) to apply the method we need an action which also admits an equivariant hyperk\"ahler moment map. Recall that if $U$ is a triholomorphic vector field then for each $i$, $\iota(U)\omega_i=d\mu^U_i$ for some function $\mu_i^U$ and putting the moment maps $\mu_i^U$ together we get the hyperk\"ahler moment  map
$${\bf \mu}:X\rightarrow {\bf R}^n\otimes {\bf R}^3$$
which, if equivariant, represents $X$ as the total space of a principal ${\bf R}^n$ bundle over an open set of ${\bf R}^n\otimes {\bf R}^3$. The Legendre transform construction reduces the hyperk\"ahler equations in this situation to finding a real-valued function $F({\bf x}_1,\dots,{\bf x}_n)$ defined on ${\bf R}^n\otimes {\bf R}^3\rightarrow {\bf R}$ and  which satisfies the 3-dimensional Laplace equation
\begin{equation}
\Delta F(c_1{\bf x}, \dots, c_n{\bf x})=0
\label{F1}
\end{equation}
  for each $(c_1,\dots,c_n)\in {\bf R}^n$. Clearly if $F$ is defined from a holomorphic function $f(z_1,\dots,z_n)$ by
  $$F({\bf x}_1, {\bf x}_2,\dots,{\bf x}_n)=Re f(u_1+iv_1,\dots,u_n+iv_n)$$
  with ${\bf x}_j=(u_j,v_j,w_j)$ then  it satisfies this equation. In fact this is essentially the only way to obtain a solution which is invariant under translation in one of the  coordinates of each ${\bf x}_j\in {\bf R}^3$.

\begin{prop} The hyperk\"ahler metric defined by (\ref{omega1}),(\ref{omegac}) is constructed by the Legendre transform method from $F= Re {\cal F}$ where ${\cal F}$ is a holomorphic prepotential.
\end{prop}
\noindent{\it Proof:} Recall that the introduction of the holomorphic function ${\cal F}$ required the choice of a symplectic basis on $V$. We shall need the same to implement the Legendre transform method, because we need an equivariant moment map. To see this note that if $U_j=-\partial/\partial y_j$, then from (\ref{omega1}),(\ref{omegac}) the hyperk\"ahler moment map is 
\begin{eqnarray*}
\mu_1&=&\frac{\partial \phi}{\partial x_j}\\
\mu_2+i\mu_3&=&i\sum \omega_{jk}(x_k+iy_k)
\end{eqnarray*}
and this (because of the $y_k$ terms) is not equivariant for the full group ${\bf R}^{2n}$ of isometries. However, if we choose a symplectic basis so that $\omega=\sum_1^n dx_i\wedge dx_{n+i}$ and take the action of ${\bf R}^n$ generated by $U_1,\dots,U_n$, we have for $1\le j\le n$ the moment map for $U_j$
$$\mu_1=\frac{\partial \phi}{\partial x_j}=\xi_j$$
$$\mu_2+i\mu_3=-i(x_{n+j}+iy_{n+j})$$
and this is equivariant since 
$$U_j\cdot y_{n+k}=-\frac{\partial y_{n+k}}{\partial y_j}=0$$ for $1\le j \le n$.
The hyperk\"ahler moment map is now:
$${\bf \mu}(x)=(\xi_1,\dots,\xi_n,-y_{n+1},\dots,-y_{2n}, x_{n+1},\dots,x_{2n})\in {\bf R}^n\otimes {\bf R}^3$$
Since ${\cal F}$ is a holomorphic function of $w_{j}=\xi_j+ix_{n+j}$, $F= Re {\cal F}$ satisfies the equation (\ref{F1}). 

To find the hyperk\"ahler metric for such an $F$, we follow \cite{HKLR} putting $$z_j=-y_{n+j}+ix_{n+j}=i(x_{n+j}+iy_{n+j})$$ for $1\le j \le n$. These are holomorphic functions with respect to the complex structure $J_1$. According to \cite{HKLR}, a K\"ahler potential for this complex structure is
$$K=F-\sum_1^n (u_k+\bar u_k)\xi_{k}$$
where for $1\le j\le n$
$$\frac{\partial F}{\partial \xi_j}=u_j+\bar u_j$$
But from (\ref{FF}) this means that $u_j+\bar u_j=x_j$ and
\begin{equation}
K=F-\sum_1^n x_k\xi_k
\label{pot}
\end{equation}
To summarize, we have $F(x_{n+1},\dots,x_{2n},\xi_1,\dots,\xi_n)$ and from the complex Lagrangian submanifold structure of $M$ we obtain from (\ref{FF}) for $1\le j\le n$,
\begin{equation}
x_j=\frac{\partial F}{\partial \xi_j},\qquad \xi_{n+j}=-\frac{\partial F}{\partial x_{n+j}}
\label{coord}
\end{equation}
From the $\Omega_1$-Lagrangian description of $M$ we also have, for $1\le j\le 2n$,
\begin{equation}
\xi_j=\frac{\partial \phi}{\partial x_j}
\label{coordx}
\end{equation}
so in particular we can write (\ref{pot}) as
$$K=F-\sum_1^n x_k\frac{\partial \phi}{\partial x_k}$$
Differentiating with respect to $x_j$ for $1\le j \le n$ we have
\begin{eqnarray*}
\frac{\partial K}{\partial x_j}&=&\sum_1^n \frac{\partial F}{\partial \xi_k}\frac{\partial \xi_k}{\partial x_j}-\frac{\partial \phi}{\partial x_j}-\sum_1^n x_k\frac{\partial^2 \phi}{\partial x_k \partial x_j}\\
&=&\sum_1^n x_k \frac{\partial^2 \phi}{\partial x_k \partial x_j}-\frac{\partial \phi}{\partial x_j}-\sum_1^n x_k \frac{\partial^2 \phi}{\partial x_k \partial x_j}\\
&=&-\frac{\partial \phi}{\partial x_j}
\end{eqnarray*}
using (\ref{coord}) and (\ref{coordx}) and the fact that $F$ is independent of $x_j$ for $1\le j \le n$.
Similarly for $n+1\le j \le 2n$,
\begin{eqnarray*}
\frac{\partial K}{\partial x_j}&=&\sum_1^n \frac{\partial F}{\partial \xi_k}\frac{\partial \xi_k}{\partial x_j}+\frac{\partial F}{\partial x_j}-\sum_1^n x_k\frac{\partial^2 \phi}{\partial x_k \partial x_j}\\
&=&\sum_1^n x_k \frac{\partial^2 \phi}{\partial x_k \partial x_j}-\xi_j-\sum_1^n x_k \frac{\partial^2 \phi}{\partial x_k \partial x_j}\\
&=&-\frac{\partial \phi}{\partial x_j}
\end{eqnarray*}
Hence $K=-\phi$ modulo an additive constant. Since $\phi$ is a  K\"ahler potential from Proposition 5, we have the same metric (taking into account a difference of sign convention).
\vskip .25cm

\vskip .25cm
\noindent {\bf Remarks}
\noindent 1. To globalise this metric  presents some choice. One could, as in \cite{F}, define it on $T^*M$. Its local structure is, however, that of a principal bundle with structure group a translation group. As such it  has no geometrically distinguished zero section. In the context of complex Lagrangian submanifolds, it  can be defined on  the space of pairs of a complex Lagrangian submanifold together with  a  line bundle of fixed Chern class over it. In this context it is defined on an open set of Markman's moduli space \cite{Mark} of Lagrangian sheaves, which is itself an integrable system.
\vskip .15cm
\noindent 2. In \cite{H} a K\"ahler metric on the moduli space of pairs consisting of a special Lagrangian submanifold of a Calabi-Yau manifold together with a flat line bundle was defined and conjectured to be itself Calabi-Yau. In the case that the Calabi-Yau is hyperk\"ahler and the special Lagrangian submanifold is complex Lagrangian with respect to one of the complex structures, this metric is precisely the one defined above. Since it is hyperk\"ahler it is {\it a fortiori} Calabi-Yau.
\vskip .15cm
\noindent 3. From the point of view put forward in this paper we have travelled from a complex Lagrangian submanifold of ${\bf C}^{2n}$ (Remark 1 of Section 4) to a hyperk\"ahler metric. This is essentially the route followed by Cortes in \cite{Cort}.

\end{document}